\documentclass{article}
\usepackage[a4paper,left=3cm,top=3.5cm]{geometry}
\usepackage{verbatim}
\usepackage{graphicx}
\usepackage{rotating}
\usepackage{tikz}
\usepackage[all,cmtip]{xy}
\usepackage{latexsym,bm,easymat}
\usepackage{amsmath,,calligra,mathrsfs,amscd,amssymb,amsfonts,mathdots,ntheorem}
\usepackage{indentfirst}
\usepackage[colorlinks,linkcolor=blue,citecolor=red]{hyperref}
\usepackage{xcolor,graphicx,blkarray}
\usepackage{latexsym}
\usepackage{longtable}
\usepackage{colortab}
\usepackage{cases}
\usepackage[titletoc]{appendix}
\begin{document}

%%%%%%%%%%%
\def\comp{\ensuremath\mathop{\scalebox{.6}{$\circ$}}}
\def\QEDclosed{\mbox{\rule[0pt]{1.3ex}{1.3ex}}} % 定义实心符
\def\QEDopen{{\setlength{\fboxsep}{0pt}\setlength{\fboxrule}{0.2pt}\fbox{\rule[0pt]{0pt}{1.3ex}\rule[0pt]{1.3ex}{0pt}}}} %定义空心符
\def\QED{\QEDopen} % 选填\QEDclosed得到实心
\def\pf{\noindent{\bf Proof}} %定义证明，注意选择中英文
\def\endpf{\hspace*{\fill}~\QED\par\endtrivlist\unskip \hfill}
\def\Hom{\mbox{Hom}}
\def\ch{\mbox{ch}}
\def\Str{\mbox{Str}}
\def\End{\mbox{End}}
\def\Re{\mbox{Re}}
\def\Im{\mbox{Im}}
\def\tri{\triangle}
\def\cs{\mbox{cs}}
\def\mcx{\mathcal{X}}

%%%%%%%%%%
\theoremstyle{mystyle}
\renewcommand{\contentsname}{\center{Content}}
\newtheorem{defi}{Definition}[section]
\newtheorem{thm}[defi]{\textbf{Theorem}}
\newtheorem{conj}[defi]{\textbf{Conjecture}}
\newtheorem{cor}[defi]{\textbf{Corollary}}
\newtheorem{lemma}[defi]{\textbf{Lemma}}
\newtheorem{exa}[defi]{\textbf{Example}}
\newtheorem{prop}[defi]{Proposition}
\newtheorem{rmk}[defi]{Remark}
\newtheorem*{ack}{Acknowledgment}
\newtheorem*{emp}{}
\newtheorem*{que}{Question}
 
\def\mcR{\mathcal{R}}

\def\pari{\partial_i}
\def\parj{\partial_j}
\def\park{\partial_k}
\def\parl{\partial_l}
\def\parp{\partial_p}
\def\parq{\partial_q}
\def\para{\partial_a}
\def\parb{\partial_b}
\def\imnum{\sqrt{-1}}

\def\nai{\nabla_i}
\def\nj{\nabla_j}
\def\np{\nabla_p}
\def\nl{\nabla_l}
\def\nk{\nabla_k}
\def\nr{\nabla_r}
\def\na{\nabla_a}
\def\nb{\nabla_b}

\def\mcP{\mathcal{P}}
\def\mcQ{\mathcal{Q}}
 %Let $(X,J)$ be  a compact almost complex manifold. By the Serre duality and the above lemma, the conjugate map $H^{2,0}_{Dol}\to H^{0,2}_{Dol}$ is an  isomorphism.\end{cor}

\def\lapg{\Delta_{g}}
\def\lapgn{\Delta_{g'} }
\def\trgng{\mbox{tr}_g(g') }

\def\parr{\partial_r}
\def\parl{\partial_l}

\def\ns{\nabla_s}

%The fifth term can be expressed as 

\def\paral{\partial_\alpha}
\def\parbet{\partial_\beta}
\def\zalp{Z_\alpha}
\def\zbet{Z_\beta}

\def\galbe{G_{\alpha\bar\beta}}
\def\bzalp{\overline{Z}_\alpha}
\def\bzbet{\overline{Z}_\beta}
\def\alpbbet{{\alpha\bar\beta}}
\def\alpbga{{\alpha\bar\gamma}}
\def\sigbbe{{\sigma\bar\beta}}
\def\sigbga{{\sigma\bar\gamma}}

\def\gabbet{\gamma\bar\beta}
\def\parga{{\partial_\gamma}}
\def\paralp{{\partial_\alpha}}
\def\parbet{{\partial_\beta}}
\def\alpbet{{\alpha\beta}}
\def\alpgam{{\alpha\gamma}}
\def\gamalp{{\gamma\alpha}}
\def\gami{{\gamma i}}
\def\beti{{\beta i}}
\def\alpi{{\alpha i}}
\def\betp{{\beta p}}
\def\gambet{{\gamma\beta}}

\def\tmd{\mathcal{\tilde{D}}^+_J}
\def\mcP{\mathcal{P}}
\def\mcPe{\mcP_\epsilon}
\def\omo{\omega_1}
\def\tom{\Tilde{\omega}_1}
\def\dj{d_J}
\def\pbp{\partial\bar\partial}
\def\tmw{{\mathcal{\widetilde W}}}

 \def\mbn{{\mu\bar\nu}}
 \def\gambdelt{{\gamma\bar\delta}}

\def\BR{{B_R(0)}}

\title{A note on $\tmd$-operator}
\author{Dexie Lin, Hengyu Zhou}
\date{}%No. 1120140012}
\maketitle

\begin{abstract}    In almost K\"ahler manifolds, one of the challenges is to construct an elliptic operator on functions that plays a role analogous to the $\partial\bar{\partial}$ operator in complex or K\"ahler manifolds. One of the aims of this paper is to revisit the $\tmd$-operator introduced in \cite{TWZZ}. We will provide some local analysis estimates and highlight several difficulties that remain to be addressed. Additionally, we use the Atiyah-Hitchin-Singer operator to demonstrate that every $d$-exact $(1,1)$-form is globally $\tmd$-exact for any compact taming symplectic $4$-manifold.
\end{abstract}

%\tableofcontents

\section{Introduction}
\def\mwd{\mathcal{W}_d}
\def\mw{\mathcal{W}}

One prominent direction in the study of symplectic manifolds involves extending equations from K\"ahler manifolds to   symplectic manifolds. Among these, the Calabi-Yau equation plays an   important role in the study of symplectic $4$-manifolds. This is programmed by Donaldson in \cite{Don06}. A major challenge in solving this equation lies in the absence of the $\partial\bar{\partial}$-lemma in the symplectic manifolds. 
\begin{itemize}
    \item Weinkove \cite[Section 2-Almost Kähler potential]{Wein07} introduced the concept of almost Kähler potentials for compact almost Kähler $4$-manifolds.
    \item Around the same time, Lejmi \cite[Theorem 2.1]{Lejmi06} observed that any almost complex $4$-manifold admits a locally compatible symplectic form. In his proof, Lejmi defined an elliptic operator on almost complex $4$-manifolds. 
\end{itemize}
 Based on these contributions, Tan, Wang, Zhou, and Zhu \cite[Definition 2.5]{TWZZ} combined the elliptic operators introduced by Weinkove and Lejmi to define a generalized $\sqrt{-1}\partial\bar{\partial}$-operator, called the $\tilde{D}^+_J$-operator, for compact tamed symplectic $4$-manifolds. In their work, they made progress toward resolving Donaldson's taming conjecture \cite{Don06} under a certain condition ($h^-_J = b^+ - 1$).
Additionally, Tan, Wang, Zhou, and Zhu introduced the $\mw$-operator, which can be interpreted as the $Jd$-operator augmented with additional terms. This naturally raises the following question:

\begin{que}
    Let $a$ be a real-valued $1$-form satisfying $d^-_J a = 0$. Does there exist a real-valued function $f$ such that $\tmd(f) = da$?
\end{que}
A classical approach to such problems is the $L^2$-method, as seen in \cite{Hor65}, analogous to the $\bar{\partial}$-problem solved by H\"ormander. During a paper-reading seminar, the authors observed certain challenges related to the local analysis of the $\tmd$-operator, as detailed in Lemma \ref{lemma-local} and Remark \ref{rmk-local-1}. One of the goals of this paper is to outline these difficulties concerning the local exactness of the $\tmd$-operator, as discussed in Lemma \ref{lemma-local}.
Furthermore, we employ the Atiyah-Hitchin-Singer operator \cite{AHS} to address the aforementioned question. Our main result is the following:

\begin{thm}\label{thm-1}
    Let $(M, J, g, F)$ be a compact almost Hermitian $4$-manifold, where $J$ is tamed by a symplectic form $\omega$. If a real-valued $(1,1)$-form $\psi$ can be expressed as $\psi = da$ for some real-valued $1$-form $a$, then $\psi$ is $\tmd$-exact.
\end{thm}
Here, a compact manifold mean a compact manifolds without boundary.

\noindent
{\bf Key Words:} taming symplectic condition, $4$-manifolds, $\tmd$-operator 
\begin{ack}
The authors would like to thank Yi Wang for the discussion.
     The first author is supported by NSFC Grant  No. 12301061.
\end{ack}
 
%\begin{ack}   This note is supported by NSFC No. 12301061.\end{ack}
\section{Preliminary on $\tilde{D}^+_J$-operator on tamed almost complex $4$ manifold}
\def\slc{sl_2(\mathbb C)}
\def\mbr{\mathbb R}
\def\mbc{\mathbb C}
\def\mfa{\mathfrak a}
\def\mcz{\mathcal{Z}}
\def\cinf{C^\infty}
\def\mcA{\mathcal{A}}
\def\mcH{\mathcal{H}}

In this section, we reintroduce the operators $\tmd$ and $\tmw$ on a compact taming symplectic $4$ manifold  {from \cite{TWZZ}}.
First, we recall two basic definitions. 
\begin{defi}
	For a smooth manifold $M$ of real dimension $2n$. If there exists a smooth section $J\in\Gamma(End(TM))$ on $TM$  satisfying $J^2=-1$, then we call
	$(M,J)$ an almost complex manifold with almost complex structure $J$.
\end{defi}
\begin{defi}
	A metric $g$   on an almost complex manifold $(M,J)$ is called Hermitian, if $g$ is $J$-invariant, i.e. $g(J-,J-)=g(-,-)$.
	The imaginary component $\omega$ of $g$ defined by $\omega=g(J-,-)$, is a real non-degenerate $(1,1)$-form. 
	The quadruple $(M,J,g,\omega)$ is called an \textit{almost Hermitian} manifold. 
	Moreover, when $\omega$ is $d$-closed, the quadruple $(M,J,g,\omega)$ is called  an
	\textit{almost  K\"ahler} manifold. 
\end{defi}
Let  $T_{\mathbb C}M=TM\otimes_{\mathbb R}\mathbb C$ be the complexified tangent bundle. We have the decomposition with respect to $J$,
\[T_{\mathbb C}M=T^{1,0}M\oplus T^{0,1}M,\]
where $T^{1,0}M$ and $T^{0,1}M$ are the eigenspaces of $J$ corresponding to the eigenvalues $i$ and $-i$ respectively.
The almost complex structure $J$ acts on the cotangent bundle $T^* M$ by $J  \alpha(v)=\alpha(J v)$, where $\alpha$ is a $1$-form and $v$ a vector field on $M$. This action $J$  can be extended to any $p$-form $\psi$ by,
\begin{eqnarray}
  (J\psi)\left(v_1, \cdots, v_p\right)= \psi\left(J v_1, \cdots, J v_p\right).  \label{formu-J-action}
\end{eqnarray}
Analogously, we decompose $T^*_{\mathbb C}M:=T^{*}M\otimes_{\mathbb R}\mathbb C$ as
$T^*_{\mathbb C}M=T^{*,(1,0)}_{\mathbb C}M\oplus T^{*,(0,1)}_{\mathbb C}M$.
%where $T'^*_{\mathbb C}M$ denote the eigenspace of $\sqrt{-1}$ and $T''^*_{\mathbb C}M$ denote the eigenspace of $-\sqrt{-1}$.
Denoting $ \mathcal A^*(M)$ by the space of sections of  $\bigwedge^*T^*_{\mathbb C}M$,   the decomposition holds:
\[ \mathcal A^r(M)=\bigoplus_{p+q=r} \mathcal A^{p,q}(M),\]
where $\mathcal{A}^{p,q}(M)$ denotes the space  of   forms of $\bigwedge^{p,q}T^*_{\mathbb C}M:=\bigwedge^pT^{*,(1,0)}_{\mathbb C}M\otimes\bigwedge^qT^{*,(0,1)}_{\mathbb C}M$.
The action \eqref{formu-J-action} gives the identity,
\[ J^2\big|_{\mathcal A^2(M)}=1.\]
In particular,
 \[ J\big|_{\mathcal A^{1,1}(M)}=1 \mbox{ and } J\big|_{\mathcal{A}^{2,0}\oplus\mathcal A^{0,2}}=-1.\]
Let $(M,J)$ be a compact almost complex  manifold. We set
$\mcz^2_\mbr$ as the space of real closed $2$-forms and  $\mathcal{A}^\pm_J$ as the $\pm1$-eigenspaces of $J$ on $\mcA^2$.  We define
\[\mcz^\pm_J=\mcA^\pm_J\cap \mcz^2_\mbr.\]

\begin{defi}[Draghici Li and Zhang \cite{DLiZh10}]
    Set 
\[H^\pm_J=\{\mfa\in H^2(M;\mathbb R)\mid \exists a\in\mcz^\pm_J,~[a]=\mfa\}.\]
$J$ is said  to be $C^\infty$-pure if $H^+_J\cap H^-_J=\{0\}$, and $C^\infty$-full if $H^+_J+H^-_J=H^2(M;\mathbb R)$.
\end{defi}

%\begin{thm} If $(M,J)$ is a compact almost complex $4$ manifold. Then, $J$ is $C^\infty$-pure and full.   \end{thm}
\noindent
%Denote by $h^\pm_J=\dim H^\pm_J$. Let $\mathcal{A}^\pm$ denote the spaces self-dual and anti-self-dual $2$-forms and \[d^\pm:\mathcal{A}^1\to\mathcal{A}^2\to\mathcal{A}^\pm,\] be the corresponding projection of $d$ to  $\mathcal{A}^\pm$.By the Hodge-decomposition, it establishes\[\mcA^\pm=\mathcal{H}^\pm_g\oplus im(d^\pm),\]where $\mcH^\pm_g$ denotes the spaces of self-dual and anti-self-dual $\Delta_d$-harmonic $2$-forms.Note that \begin{equation}    d^\pm d^*:\mcA^\pm\to\mcA^\pm,\label{operator-ddpm}\end{equation}are  elliptic operators whose kernels are $\ker(d^\pm d^*)=\mcH^\pm_g$.If $d^\pm u$ is $d$-closed, then \begin{equation}   0=\int_Mdd^\pm_gu\wedge u=-\int_Md^\pm_g u\wedge d u=\mp\int_M|d^\pm u|^2,\label{eqn-d+-d-exact}\end{equation}hence $u\in\ker(d^\pm_g)=\ker(d)$.We define \[H^\pm_g=\{\mfa\in H^2(M;\mathbb R)\mid \exists a\in\mcz^\pm_J,~[a]=\mfa\}.\]
Using the Hodge decomposition, one has the following.
\begin{thm}[Draghici Li and Zhang {\cite[Theorem 2.2]{DLiZh10}}]    If $(M, J)$ is a compact almost complex $4$-manifold, then $J$ is $C^\infty$-pure and full and there is a    direct sum cohomology decomposition
    \[H^2(M;\mathbb R)=H^+_J\oplus H^-_J, \]
   i.e., 
    $b_2=h^+_J+h^-_J$. 
    In particular, the cup product of $H^2_{}(M;\mathbb R)$ is non-degenerate of type $(b^+-h^-_J,b^-)$.
\end{thm}

Note that $b^+$ and $b^-$ are topological numbers   and $b_2=b^++b^-$, however $h^\pm_J$ depends on  {the almost complex structure} $J$.
\begin{thm}[Lejmi {\cite[Theorem 2.1]{Lejmi06}}]\label{lemma-2.3}
    Let $(M,J,g,F)$ be a compact almost Hermitian $4$ manifold. We set 
    \[P:\mcA^-_J\to \mcA^-_J,~\psi\mapsto P^-_J(dd^*\psi),\]
    where $P^-_J:\mcA^2\to \mcA^-_J$ is the projection. Then, $P$ is a self-adjoint elliptic, whose kernel $\ker(P)=\mcz^-_J$, i.e.
    \[\mcA^-_J=\ker(P)\oplus im(d^-_J).\]
\end{thm}
\begin{rmk}
The operator $d^-_Jd^*$ is not elliptical for higher dimensional manifolds.  For example, on a compact K\"ahler manifold $(M,J,\omega,g)$  with complex dimension $n\geq3$. Then, one has $d^-_Jd^*\alpha=\partial\partial^*\alpha^{2,0}+\bar\partial\bar\partial^*\alpha^{0,2}$. We choose $\alpha=\partial^*\beta^{3,0}+\bar\partial^*\beta^{0,3}$, for any $\beta^{3,0}\in\mathcal{A}^{3,0}$. It is clear that such 2-forms belong to the kernel of $d^-_Jd^*$, i.e. $\ker(d^-_Jd^*)$ is of infinite dimension. 
\end{rmk} %This means $d^-Jd^*$ is not elliptic form $\mathcal{A}^-_J$ to $\mathcal{A}^-_J$, otherwise its kernel is of finite dimension. 

\begin{defi}
    Let $(M,J,g,F)$ be a compact almost Hermitian $4$-manifold. $J$ is tamed by a symplectic form $\omega$, if the $(1,1)$-component of $\omega$ equals to $F$, that is,  
    $\omega^{1,1}=F$.
\end{defi}

We want to find an expression of a $(1,1)$-exact form in terms of functions. 
First,
it is convenient to set the space $$L^2_2(M)_0:=\{f\in L^2_2(M)\mid \int_Mfdvol=0\}.$$
For $f\in L^2_2(M)_0$, it is clear $d^-_JJdf\perp \ker(P)$.  Since $P=d^-_Jd^*$ is a formally self-adjoint elliptic operator, there exists $\sigma^1_f\in L^2_2(\mathcal{A}^-_J)$ such that 
\begin{equation}
d^-_JJdf+d^-_Jd^*(\sigma^1_f)=0.\label{eqn-2.26}
\end{equation}
%\[\]
On the other hand, for any $d$-exact $\psi\in L^2(\mathcal A^{1,1}_\mbr(M))$, there exists $u_\psi\in L^2_1(\mathcal{A}^{1}_\mbr)$ such that 
\[\psi=d(u_\psi )=d^+_J(u_\psi ),~i.e., d^-_J(u_\psi )=0.\]
Hence, equation \eqref{eqn-2.26} yields 
\begin{equation}
d^+_J(u_f )=d^+_J(Jdf+ d^*(\sigma^1_f)),\label{eqn-d+J}
\end{equation} 
for some $u_f\in L^2_1(\mathcal{A}^{1}_\mbr)$.
Next,
consider the self-adjoint elliptic operator 
\[d^+d^*=\frac12\Delta_d|_{\mathcal{A}^+_\mbr}:\mathcal{A}^+_\mbr\to\mcA^+_\mbr,\]
where $d^+=\frac{1}{2}(d+*d)$.
There exists $\sigma^2_f\in\mathcal{A}^{-}_J$ satisfying
\begin{equation*}
d^+(u_f )=d^+d^*(f\omega_1+\sigma^1_f+\sigma^2_f).  
\end{equation*}
It computes 
\begin{eqnarray*}
d^*(f\omega_1)&=&-*d(f\omega_1)\\
&=&-*(df\wedge\omega_1)\\
&=&-*(df\wedge F+df\wedge d^-_J(v ))\\
&=&Jdf-*(df\wedge d^-_J(v )),
\end{eqnarray*}
where we use the identity $*J\alpha=\alpha\wedge F$ for any $1$-form $\alpha$ on the last equality. 
By the above calculation, it holds that 
\begin{eqnarray*}
\psi_f&:=&d(u_f )\\
&=&dd^*(f\omega_1+\sigma^1_f+\sigma^2_f)\\
&=&dJdf+dd^*(\sigma^1_f)-d*(df\wedge d^-_J(v ))+dd^*(\sigma^2_f).
\end{eqnarray*}
Combining the equations \eqref{eqn-d+J} and \eqref{eqn-2.26}, we write 
\begin{equation}
\begin{cases}
	d^-_JJdf+d^-_Jd^*(\sigma^1_f)=0,\\
	d^-_J*(df\wedge d^-_J(v ))-d^-_Jd^*(\sigma^1_f)=0.
\end{cases}\label{eqn-2.29}
\end{equation}
Summarizing the above arguments, one has the following operators on a compact tamed symplectic $4$-manifold.
\begin{defi}[Tan Wang Zhou and Zhu {\cite[Definition 2.5]{TWZZ}}]
    Set 
    $$\mathcal{W}:L^2_2(M)_0\to L^2_1(\mathcal{A}^1_\mbr),
    ~f\mapsto Jdf+d^*(\sigma^1_f).$$ 
    It satisfies  $d^-_J\mathcal{W}(f)=0$. Define $D^+_J:L^2_2(M)_0\to L^2(\mcA^{1,1}),~f\mapsto d\mathcal{W}(f)$.
    Set 
    \[
   \mathcal{ \Tilde  W}:L^2_2(M)_0\to L^2_1(\mathcal{A}^1_\mbr),~f\mapsto \mathcal{W}(f)-*(df\wedge d^-_J(v ))+d^*(\sigma^2_f).
    \]
    It satisfies $d^*\tilde{\mathcal{W}}(f)=0,~d^-_J\tilde{\mathcal{W}}(f)=0$. 
    Define $$\tilde{D}^+_J:L^2_2(M)_0\to L^2(\mcA^{1,1}_\mbr),~f\mapsto d\Tilde{\mathcal{W}}(f).$$
\end{defi}
Observe that for compact K\"ahler surface, we choose 
$\mw f=\tmw f=Jdf$ and $\tmd=2\imnum\partial\bar\partial$.

It is convenient to consider the formal adjoint of $\tmw$.
For any $a\in\mathcal{A}^1_\mbr$ satisfying $d^*a=0=d^-_Ja$, it derives 
\begin{eqnarray*}
	(\tmw(f),a)_{L^2}&=&-\int_Ma\wedge d(f\omega_1+(\sigma^1_f)+(\sigma^2_f))\\
	&=&-\int_Mda\wedge(f\omega_1+(\sigma^1_f)+(\sigma^2_f))\\
	&=&-\int_Md^+_J(a)\wedge fF=(f,\tmw^*(a))_{L^2}.
\end{eqnarray*}
That is
\begin{eqnarray}
    {\tmw}^*(a)=-\Lambda d^+_Ja,\label{formu-tW-adjoint}
\end{eqnarray}
where $\Lambda$ denotes the contraction with the form $F$. 
Analogously, the calculation for $d^-_Ja=0$,
\begin{eqnarray*}
    (\mathcal{W}(f), a) & =&\left(-*(d f \wedge F)+d^*(\sigma_f^1), a\right)) \\
& =&\left(-*(d(f F)-f d F)+d^*\left(\sigma_f^1\right), a\right) \\
&  =&\left(-*(d(f F)-f d F) ), a\right) +(\sigma^1_f,da)\\
& =&\left(d^*(f F), a\right)+\left(* (f d F), a\right)+ (\sigma^1_f,d^+_Ja) \\
& =&((f F), d a)+\left(* (f d F), a\right) \\
& =&\int_M f F \wedge d a-f d F \wedge a ,
\end{eqnarray*}
gives the formally adjoint operator
\begin{equation}
    \mathcal{W}^*(a)=\frac{2\left(d_J^{+} A \wedge F+a \wedge d F\right)}{F^2}\label{formu-W-adjoint}
\end{equation}

One also has the following. 

\begin{lemma}[Tan Wamg Zhou and Zhu {\cite[Lemma 3.1 and 3.2]{TWZZ}}]
    On a compact almost Hermitian $4$ manifold,
$d^+_J:L^2_1(\mcA^1_\mbr)\to L^2(\mcA^{1,1})$ has a closed range and  $\Tilde{D}^+_J:L^2_2(M)_0\to L^2(\mathcal{A}^{1,1})$ has a closed range. 
\end{lemma}

\section{Local behavior on $\tmd$-operator}
In this section, we give some difficulties which are crucial in the proof of taming conjecture, cf. \cite[Proof of Theorem 1.1 and Proof of Appendix Theorem C.12]{TWZZ}. 
\subsection{Local exactness}
In this subsection, we focus on the local property of $\tmw$ and $\tmd$ in \cite[Appendix A.3]{TWZZ}.
We consider the local model on $\tmd$ and $\tmw$. In \cite[Appendix A.3]{TWZZ}, Tang Wang Zhou and Zhu show the following short exact complex
$$
L^2(\Omega)_0 \xrightarrow{\tmw} L^2(\mcA^1_{\mathbb{R}}  (\Omega)) \xrightarrow{d_J^{-}}  L^2(\mcA^-_J(\Omega)),
$$
where $\Omega$ is a $J$-pseudoconvex domain  in the sense of \cite[Corollary 6.5]{HV15} or \cite[P. 930]{HLP16} and 
$L^2(\Omega)_0=\{f\in L^2(\Omega)\mid \int_\Omega f dvol=0\}$.
This is equivalent to  whether the kernel of $d_J^{-}$ is equal to the image of $\widetilde{\mathcal{W}}$. We apply  the method of $L^2$ estimates for the complex  $\bar{\partial}$-problem(cf. \cite[Chapter IV]{Hor90}) to the approach. Such a approach is also given in \cite[Appendix A]{TWZZ}.
 For simplification, we summarize the above discussion in terms of the model of Hilbert spaces as following: 
$$
H_1 \xrightarrow{T} H_2 \xrightarrow{S} H_3,
$$
where $H_1, H_2, H_3$ are all $L^2$ spaces, and $T, S$ are linear, closed, and densely defined closed operators. Assume   $S\comp T=0$, the problem is, whether for any $  g \in \operatorname{ker} S$,  there is a solution to
$$
T f=g.
$$
 First, notice a simple fact that $T f=g$ is equivalent to$$(T f, h)_{H_2}=(g, h)_{H_2}, \forall h \in \text { D },$$
where $D$ is a dense subset of $H_2$. 
This is because that $(T f-g, h)_{H_2}=0, \forall h \in D$ is equivalent to  $\left(T f-g, H_2\right)_{H_2}=0,\mbox{ i.e., }T f=g$
.Let $T^*$ be an adjoint operator of $T$ in the sense of distributions. By functional analysis theory, $T^*$ is a closed operator and $\left(T^*\right)^*=T$ if and only if $T$ is closed.
One has a well-known model.
\begin{thm}[H\"ormander {\cite[Theorem 1.1.4]{Hor65}}]\label{lemma-partial-bar}
Let \[ H_1\overset{T}{\longrightarrow}H_2\overset{S}{\longrightarrow}H_3 \]
be a complex ($ST=0$) over Hilbert spaces $H_1,H_2,H_3$, where $T$ and $S$ are densely defined operators with closed graphs, whose domains are denoted by $D_T$ and $D_S$ respectively. Let $V$ be a closed subspace of $H_2$
containing $im(T)$.
Suppose that the estimate
\[\|h\|_{H_2}\leq C(\|T^*h\|_{H_1}+\|Sh\|_{H_3}),\] holds
for any $h\in D_{T^*}\cap D_S\cap V$, where $D_{T^*}$ is the domain of $T^*$ and $C$ is a constant. Then, the following hold:
\begin{itemize}
    \item[(1)]   For any $v \in \ker S \cap V$ one can find $w \in D_T$ such that $T w=v$ and $\|w\|_{H_1} \leq C\|v\|_{H_2}$.
 \item[(2)] For any $w \in imT^*$ one can find $v \in D_{T^*}$ such that $T^* v=w$ and $\|v\|_{H_2} \leq$ $C\|w\|_{H_1}$. 
\end{itemize}
%$|(g,h)_{H_2}|\leq C\|g\|_{H_2}\|T^*h\|_{H_1}$ for any $g\in \ker(S)\cap V$ and $h\in D_{T^*}\cap D_S$, and the equation $T(f)=g$ has a solution $f\in D_T$ to $g\in \ker(S)$. Moreover, one can find a   solution  $f\in\ker(T)^\perp$ with the estimate $\|f\|_{H_1}\leq C\|g\|_{H_2}$.
\end{thm}
%From (A.11),
%$(T f, h)_{H_2}=(g, h)_{H_2}, \forall h \in$ some dense subset. If this dense subset is contained in $D_{T^*}$, then, noticing $(T f, h)_{H_2}=\left(f, T^* h\right)_{H_1}$,

%The formula \eqref{eqn-w-adjoint} implies 
Turn back to our case. 
Let $\varphi$ be a strictly $J$-plurisubharmonic function(\cite[Definition 3.1]{HV15}), i.e. 
\[\sqrt{-1}\partial\bar\partial \varphi>0,~\mbox{ in }\Omega.\]
Define 
\[ H_1=L^2_{ \varphi}(\Omega):=\{f\mid \int_\Omega |f|^2e^{-\varphi}<\infty\}. \]
With the similar notations, we define 
\[H_2=  L^2_{\varphi}(\mcA^1_{\mathbb{R}}  (\Omega)) \mbox{ and } H_3= L^2_{\varphi}(\mcA^-_J(\Omega).\]
Parallel to the calculation of \eqref{formu-tW-adjoint} 
we  directly obtain
$$
\widetilde{\mathcal{W}}^*_{\varphi} a=\frac{-2 F \wedge d_J^{+}\left(e^{-\varphi} a\right)}{F^2} \cdot e^{\varphi} .
$$
We need to calculate   $\left\|\widetilde{\mathcal{W}}_{\varphi}^* a\right\|_{H_1}^2+\left\|d_J^{-} a\right\|_{H_3}^2$, for $a \in D_{\widetilde{\mathcal{W}}_{\varphi}^*} \cap D_{d_J^{-}} \cap \mathcal A_{\mathbb{R}}^1(\bar{\Omega})$, where $ D_{\widetilde{\mathcal{W}}_{\varphi}^*}$ and $D_{d_J^{-}}$ denote the domains of ${\widetilde{\mathcal{W}}_{\varphi}^*}$ and ${d_J^{-}}$ respectively. 
Choose a local   unitary frame $\left\{e_1, e_2\right\}$ for $T^{1,0}(\Omega)$ and a local complex coordinate $\left\{z^1, z^2\right\}$ in a neighborhood of $p$ satisfying $e_i(p)=\left.\frac{\partial}{\partial z^i}\right|_p$ with respect to the Hermitian metric. Denote $\left\{\theta^1, \theta^2\right\}$ by the dual frame of $\left\{e_1, e_2\right\}$. Hence, one has 
%$$h=g_J-\sqrt{-1} F=\theta_1 \otimes \bar{\theta}_1+\theta_2 \otimes \bar{\theta}_2$$and
$$
F=\sqrt{-1}(\theta^1 \wedge \bar{\theta}^1+\theta^2 \wedge \bar{\theta}^2 ).
$$
We write the structure equation as 
$$\bar\mu \theta^s= \sum_{k<j}N^s_{\bar k\bar j}\bar\theta^k\wedge\bar\theta^j,~\bar\partial \theta^s=-\sum_{k,j}C^{s}_{k\bar j}\theta^k\wedge \bar\theta^j,~
\partial \theta^s=- \sum_{k<j}C^s_{kj}\theta^k\wedge\theta^j.$$
The coefficients of   Nijenhuis tensor $N$ is   given by
$$
N\left(e_{\bar{j}}, e_{\bar{k}}\right)=-\left[e_{\bar{j}}, e_k\right]^{(1,0)}= \sum_t N_{j k}^t e_t,~ N\left(e_j, e_k\right)=-\left[e_j, e_k\right]^{(0,1)}=\sum_t\overline{N_{jk}^t} e_{\bar{t}},
$$
and   the structure coefficients of the Lie bracket is expressed as $$[e_j,\bar e_{k}]=\sum_k C^r_{j\bar k}e_r+
\sum_lC^{\bar l}_{i\bar j}e_{\bar l}.$$
%In their paper \cite[Appendix A.3]{TWZZ}, they use
Under these notations, we write 
\begin{eqnarray}
\partial\bar\partial \phi=\sum_{i,j}(\bar e_{j}e_i \phi+[e_i,\bar e_{j}]^{1,0}\phi)\theta^i\wedge\bar\theta^j=
\sum_{i,j}(e_i\bar e_j\phi-\sum_l C^{\bar l}_{i\bar j}\cdot\bar e_l\phi)\theta^i\wedge\bar\theta^j,\label{eqn-ddc}
\end{eqnarray}
for any function $\phi$.

For $u=u_1 \bar{\theta}^1+u_2 \bar{\theta}^2$, a direct calculation yields
\begin{eqnarray*}
\partial\left(e^{-\varphi} u\right)\wedge F
&=&-e^{-\varphi}(e_1\varphi\cdot \bar\theta^1+e_2\varphi\cdot\bar\theta^2)\wedge u\wedge F+
e^{-\varphi}\partial(u_1\theta^1+u_2\theta^2)\wedge F\\
&=&-e^{-\varphi}(e_1\varphi\cdot \theta^1+e_2\varphi\cdot\theta^2)\wedge u\wedge F+
e^{-\varphi}(e_1u_1 \theta^1\wedge\bar\theta^1+e_2u_2\theta^2\wedge\bar\theta^2)\wedge F\\
&&-e^{-\varphi}(u_1 C^{\bar 1}_{1\bar1}\theta^1\wedge\bar\theta^1+u_1 C^{\bar1}_{2\bar 2}\theta^2\wedge\bar\theta^2)\wedge F-e^{-\varphi}(u_2C^{\bar 2}_{1\bar 1}\theta^1\wedge\bar\theta^1+u_2  C^{\bar 2}_{2\bar2}\theta^2\wedge \bar\theta^2)\wedge F\\
&=&-\frac{1}{2}e^{-\varphi}(e_1\varphi u_1+e_2\varphi u_2)F^2+\frac{1}{2}e^{-\varphi}(e_1u_1+e_2u_2)
F^2\\
&&-\frac{1}{2}e^{-\varphi}(u_1C^{\bar 1}_{1\bar1}+u_1 C^{\bar1}_{2\bar 2}+u_2C^{\bar 2}_{1\bar 1}+u_2 C^{\bar2}_{2\bar2})F^2.
\end{eqnarray*}
Setting  $a=u+\bar{u}$, one  writes 
\begin{eqnarray*}
d_J^{+}\left(e^{-\varphi} a\right) \wedge F & =&\left[\partial\left(e^{-\varphi} u\right)+\bar{\partial}\left(e^{-\varphi} \bar{u}\right)\right] \wedge F \\
&=&-\frac{1}{2}e^{-\varphi}(e_1\varphi u_1+e_2\varphi u_2)F^2+\frac{1}{2}e^{-\varphi}(e_1u_1+e_2u_2)
F^2\\
&&-\frac{1}{2}e^{-\varphi}(u_1C^{\bar 1}_{1\bar1}+u_1C^{\bar 1}_{2\bar 2}+u_2C^{\bar 2}_{1\bar 1}+u_2C^{\bar 2}_{2\bar2})F^2\\
&&-\frac{1}{2}e^{-\varphi}(\bar e_1\varphi \bar u_1+\bar e_2\varphi \bar u_2)F^2+\frac{1}{2}e^{-\varphi}(\bar e_1 \bar u_1+\bar e_2\bar u_2)
F^2\\
&&-\frac{1}{2}e^{-\varphi}(\bar u_1 C^{1}_{1\bar1}+\bar u_1 C^{1}_{2\bar 2}+\bar u_2 C^2_{1\bar 1}+\bar u_2C^2_{2\bar2})F^2.
\end{eqnarray*}
Summarizing the above calculation, we obtain
\begin{eqnarray*}
\widetilde{\mathcal{W}}^* a&=&
e_1\varphi u_1+e_2\varphi u_2-e_1u_1-e_2u_2+\bar e_1\varphi \bar u_1+\bar e_2\varphi \bar u_2-\bar e_1\bar u_1-\bar e_2\bar u_2\\
&&+(u_1C^{\bar 1}_{1\bar1}+u_1C^{\bar 1}_{2\bar 2}+u_2 C^{\bar2}_{1\bar 1}+u_2\bar C^2_{2\bar2})+(\bar u_1 C^1_{1\bar1}+\bar u_1 C^{1}_{2\bar 2}+\bar u_2 C^2_{1\bar 1}+\bar u_2C^2_{2\bar2}),
\end{eqnarray*}
and 
\begin{eqnarray*}
\|\widetilde{\mathcal{W}}^* a\|^2=\int_{\Omega}|\sum_j\delta_j u_j|^2e^{-\varphi} ,
\end{eqnarray*}
where $\delta_ju_j=e_ju_j-e_j\varphi\cdot u_j+\sum_k C^{\bar j}_{k\bar k}\cdot u_j$
%$$\widetilde{\mathcal{W}}^* A=\frac{\partial \varphi}{\partial z^1} u^1+\frac{\partial \varphi}{\partial z^2} u^2-\frac{\partial u^1}{\partial z^1}-\frac{\partial u^2}{\partial z^2}+\frac{\partial \varphi}{\partial \bar{z}^1} \bar{u}^1+\frac{\partial \varphi}{\partial \bar{z}^2} \bar{u}^2-\frac{\partial \bar{u}^1}{\partial \bar{z}^1}-\frac{\partial \bar{u}^2}{\partial \bar{z}^2} .$$

We need the following lemma to assist the subsequence computation. 
\begin{lemma}\label{lemma-div(L)}
If the boundary $\partial \Omega=\{r=0\}$ of a bounded domain $\Omega=\{r<0\}$ is differentiable, $|d r|=1$, and $L=\sum a_i \frac{\partial}{\partial x_i}$ is a differentiable operator of first-order, then
$$
\int_{\Omega} L f=\int_{\partial \Omega}(L r) f-\int_{\Omega}div(L)f
$$
where $div(L)=\sum_j\frac{\partial a_j(x)}{\partial x_j}$ and $(x^1,...x^n)$ is the fixed coordianate in $\Omega$. %\textcolor{red}{In your computation, $\Omega$ is a Pesudo domain in the standard complex space $\mathbb C^n$ or $\mathbb R^n$? The divergence should depend on the metric. }
\end{lemma}
%When the coefficients of $L$ are constant, the proof is using the Stokes theorem(c.f. \cite[Proposition 3.4]{CSZ}). In fact, the proof is no essentially different to \cite[Proposition 3.4]{CSZ}.

\begin{pf}
Using the integral by part, one writes 
\begin{eqnarray*}
\int_{\Omega}Lf=\int_{\Omega}\sum_j\frac{\partial}{\partial x_j} (a_j(x)f)-\int_{\Omega}\sum_j\frac{\partial a_j(x)}{\partial x_j} f\\
=\int_{\partial \Omega}(L r) f-\int_{\Omega}div(L)f,
\end{eqnarray*}
where we apply the Lemma \ref{lemma-L} below to the second equality. 
\end{pf}
\begin{lemma}\label{lemma-L}
  If the boundary $\partial \Omega=\{r=0\}$ of a bounded domain $\Omega=\{r<0\}$ is differentiable, $|d r|=1$, and $L=\sum a_i \frac{\partial}{\partial x_i}$ is a first-order differentiable operator with constant coefficients, then
$$
\int_{\Omega} L f=\int_{\partial \Omega}(L r) f.
$$
\end{lemma}
The proof of Lemma \ref{lemma-L} is to use the Stokes-lemma(cf. \cite[Proposition 3.4]{CSZ}). In fact, 
the proof of Lemma \ref{lemma-div(L)} is no essentially different to \cite[Proposition 3.4]{CSZ}. %and restated in \cite[Appendix-Proposition A.27]{TWZZ}.
\begin{rmk}
    In \cite[(A.23)]{TWZZ}, it is hard to understand the reason for the disappearance of $div(L)$.
\end{rmk}

\begin{lemma}\label{lemma-local}

    Under the  above notations, we get the following estimates:
    \begin{eqnarray}
        \left\|\widetilde{\mathcal{W}}_{\varphi}^* a\right\|_{ }^2+\left\|d_J^{-} a\right\|_{ }^2
        &=& \int_{\Omega}\sum_{i<j}|(\bar e_iu_j-\bar e_ju_i-(\bar N^i_{\bar1\bar2}u_i-\bar N^j_{\bar 1\bar 2}u_j)|^2e^{-\varphi}+\int_{i,j}\int_{\Omega}|\delta_i u_j|^2e^{\varphi}\nonumber\\
        &\geq&\sum_{i,j}\int_\Omega|\bar e_iu_j|^2e^{-\varphi}+\int_\Omega((\delta_i\bar e_j-\bar e_j\delta_i)u_i,u_j)e^{-\varphi}\nonumber\\
        &&+\sum_{i,j}\int_{\partial\Omega}(\bar e_ir\cdot(\delta_i u_j)\bar u_i-(\bar e_ir)\bar u_j\cdot(\bar e_j u_i))e^{-\varphi}\nonumber\\
        &&+\sum_{i,j}\int_{\Omega} ( div(\bar e_i)u_j,\bar e_j u_i)-(div(\bar e_j)\delta_iu_i,u_j).
        \label{esti-w-a-1}
        \end{eqnarray}
        %where we use \cite[Proposition 3.4]{CSZ} to get the first term. 
\end{lemma}

\begin{pf}
Lemma \ref{lemma-div(L)} above implies 
\begin{eqnarray*}
\int_{\Omega}f\overline{\bar e_i(u_ie^{-\varphi})}&=&-\int_{\Omega}\bar e_if\cdot \bar u_ie^{-\varphi}+
\bar e_i(f\bar u_ie^{-\varphi})\\
&=&-\int_{\Omega}\bar e_if\bar u_ie^{-\varphi}+\int_{\partial\Omega} \bar e_ir(f\bar u_i)e^{-\varphi}-\int_{\Omega}div(\bar e_i)f\bar u_ie^{-\varphi}.
\end{eqnarray*}
We  reduce the above formula to the following one:
\begin{eqnarray*}
(f_1,\delta_j f_2)_{L^2(\varphi)}&:=&\int_{\Omega} f_1\left(\overline{ e_jf_2-e_j\varphi f_2+\sum_k\bar C^j_{k\bar k} f_2}\right)e^{-\varphi}\\
&=&-\int_{\Omega}\bar e_jf_1\bar f_2e^{-\varphi}+\int_{\partial\Omega}\bar e_jrf_1\bar f_2e^{-\varphi}\\
&&-\int_{\Omega}div(\bar e_j)f_1\bar f_2e^{-\varphi}+\int_{\Omega}\sum_k C^j_{k\bar k}f_1\bar f_2e^{-\varphi}\\
&=&-(\bar e_j f_1,f_2)_{L^2(\varphi)}+(\bar e_jrf_1,f_2)_{L^2_{\partial\Omega}(\varphi)}\\
&& -(div(\bar e_j)f_1,f_2)_{L^2(\varphi)}+((\sum_k C^j_{k\bar k})f_1,f_2)_{L^2(\varphi)}.
% -(\bar e_jf_1,f_2)_{L^2(\varphi)}+((\bar e_i r)f_1,f_2)_{\partial\Omega}\\   -(div(\bar e_i)f_1,f_2)_{L^2(\varphi)}+(f_1,\sum_k\bar C^{i}_{k\bar k })_{(\varphi)}
\end{eqnarray*}
%\[(f_1,\delta_i f_2)_{L^2(\varphi)}=-(\bar e_if_1,f_2)_{L^2(\varphi)}+((\bar e_i r)f_1,f_2)_{\partial\Omega}-(div(\bar e_i)f_1,f_2)_{L^2(\varphi)}+(f_1,\sum_k\bar C^{i}_{k\bar k })_{(\varphi)},\]
for any smooth functions $f_1,f_2$. The above calculation gives
\[
\int_{\Omega}(\bar e_iu_j)\overline{(\bar e_j u_i)}e^{-\varphi}=
-(u_j,\delta_i\bar e_ju_j)_{L^2(\varphi)}+((\bar e_ir)u_j,\bar e_j u_i)_{\partial\Omega}-( div(\bar e_i)u_j,\bar e_j u_i)_{L^2(\varphi)}+(\sum_k C^i_{k\bar k}u_j,\bar e_j u_i)_{L^2(\varphi)},
\]
and 
\[
\int_{\Omega}(\delta_i u_i)\overline{\delta_j u_j}e^{-\varphi}=
-(\bar e_j\delta_i u_i,u_j)_{L^2(\varphi)}+((\bar e_jr)\delta_i u_i,u_j)_{\partial\Omega}-
(div(\bar e_j)\delta_iu_i,u_j)_{L^2(\varphi)}+(\sum_kC^j_{k\bar k}\delta_iu_i,u_j)_{L^2(\varphi)}.
\]
Expand the formula
$$
\begin{aligned}
d_J^{-}(a)= & d_J^{-}(u+\bar{u}) \\
= & \bar{\partial} u+\mu u+\partial \bar{u}+\bar\mu\bar{u} \\
= & (e_1\bar u_2-e_2\bar u_1) \theta_1 \wedge \theta_2+\left(\bar e_1u_2-\bar e_2u_1\right) \bar{\theta}_1 \wedge \bar{\theta}_2 \\
& +\left( N^2_{\bar1\bar2} \bar{u}_2-N^1_{\bar1\bar2} \bar{u}_1\right) \theta_1 \wedge \theta_2+\left(\bar N^2_{\bar 1\bar 2} u_2-\bar{N}^1_{\bar 1\bar 2} u_1\right) \bar{\theta}_1 \wedge \bar{\theta}_2.
\end{aligned}
$$It
induces 
\begin{eqnarray*}
\left\|d_J^{-} a\right\|_{ }^2&=&
\int_{\Omega}\sum_{i<j}|(\bar e_iu_j-\bar e_ju_i-(\bar N^i_{\bar1\bar2}u_i-\bar N^j_{\bar 1\bar 2}u_j))|^2e^{-\varphi}\\
&\geq&\frac{1}{2}\left(\int_\Omega(\sum_{i<j}|\bar e_iu_j-\bar e_ju_i|^2+|\bar N^i_{\bar1\bar2}u_i-\Bar{N}^j_{\bar1\bar2}u_j|^2)e^{-\varphi}\right)\\
&=&\frac{1}{2}\sum_{i,j}\int_\Omega(|\bar e_i u_j|^2-(\bar e_i u_j)\cdot(e_j\bar u_i))e^{-\varphi}+\frac{1}{2}\int_\Omega |\bar N^i_{\bar1\bar2}u_i-\Bar{N}^j_{\bar1\bar2}u_j|^2e^{-\varphi}.
\end{eqnarray*}
%where $A_{J i}$ are the coefficients of $A_J$ which is the linear operator defined in Section 2. So
%$$\begin{aligned}\left\|d_J^{-} A\right\|_{H_3}^2 & =\int_{\Omega} \sum_{i<j}\left(\left|\frac{\partial u^j}{\partial \bar{z}^i}-\frac{\partial u^i}{\partial \bar{z}^j}\right|^2+\left|A_{J j} \bar{u}^j-A_{J i} \bar{u}^i\right|^2\right) e^{-\varphi} \\& =\sum_{i, j} \int_{\Omega}\left(\left|\frac{\partial u^j}{\partial \bar{z}^i}\right|^2-\frac{\partial u^j}{\partial \bar{z}^i} \frac{\partial \bar{u}^i}{\partial z^j}\right) e^{-\varphi}+\int_{\Omega} \sum_{i<j}\left|A_{J j} \bar{u}^j-A_{J i} \bar{u}^i\right|^2 e^{-\varphi} .\end{aligned}$$

Summarizing the above results, we derive 
\begin{eqnarray*}
\left\|\widetilde{\mathcal{W}}_{\varphi}^* a\right\|_{ }^2+\left\|d_J^{-} a\right\|_{ }^2
&=& \int_{\Omega}\sum_{i<j}|(\bar e_iu_j-\bar e_ju_i-(\bar N^i_{\bar1\bar2}u_i-\bar N^j_{\bar 1\bar 2}u_j)|^2e^{-\varphi}+\int_{i,j}\int_{\Omega}|\delta_i u_j|^2e^{\varphi}\\
&\geq&\sum_{i,j}\int_\Omega|\bar e_iu_j|^2e^{-\varphi}+\int_\Omega((\delta_i\bar e_j-\bar e_j\delta_i)u_i,u_j)e^{-\varphi}\\
&&+\sum_{i,j}\int_{\partial\Omega}(\bar e_ir\cdot(\delta_i u_j)\bar u_i-(\bar e_ir)\bar u_j\cdot(\bar e_j u_i))e^{-\varphi}\\
&&+\sum_{i,j}\int_{\Omega} ( div(\bar e_i)u_j,\bar e_j u_i)-(div(\bar e_j)\delta_iu_i,u_j)
\end{eqnarray*}
where we use the identity 
\begin{eqnarray}
    (\delta_i\bar e_j-\bar e_j\delta_i)\varphi&=&e_i\bar e_j\varphi-\bar e_je_i\varphi+\bar e_je_i\varphi+\sum_{l}(\bar e_jC^{\bar i}_{l\bar l})\nonumber\\
&=&e_i\bar e_j\varphi-\sum_{r} C^{\bar r}_{i\bar j}(\bar e_{ r}\varphi)\nonumber\\
&&+\left(\sum_{r}C^r_{i\bar j}e_r+C^{\bar r}_{i\bar j}\bar e_{r} +\sum_{l}(\bar e_j C^{\bar i}_{l\bar l})\right)\varphi.\label{eqn-extra-1}
\end{eqnarray}
\[
\]
\end{pf}

\begin{rmk}\label{rmk-local-1}
   By the above calculation, the authors can hardly understand why \cite[Appendix (A.26)]{TWZZ} holds. 
   %In \cite[(A.26)]{TWZZ}, they used the notation $\partial_j=e_j$, but $\bar\partial_j\partial_i\neq \partial_i\bar\partial_j$ on almost complex case.  
The last two terms of \eqref{esti-w-a-1} are related to some first-order differential operators,  such as
\[\sum_{r}C^r_{i\bar j}e_r+C^{\bar r}_{i\bar j}\bar e_{r},\]
which are unbounded in the $L^2$-space. Moreover, the last term of the second equality in \eqref{eqn-extra-1} are also unbounded even though the first two terms are the local expression of $\partial\bar\partial \varphi$ by \eqref{eqn-ddc}, which is bounded by the hypothesis of (strict)$J$-plurisubharmonic $\varphi$.
This generates difficulty in understanding \cite[Appendix Theorem A.31 ]{TWZZ}. 
\end{rmk}

It is not hard to find a typical example with non-trivial structure coefficients. 
\begin{exa}
    Consider the local almost complex structure of $(\Omega,J)$ based on the Kodaira-Thurston $4$-manifold.
    Choose an almost complex structure given by
	\[J\frac{\partial}{\partial t}=\frac{\partial}{\partial x},~J\frac{\partial}{\partial x}=-\frac{\partial}{\partial t},~J(\frac{\partial}{\partial y}+x\frac{\partial}{\partial z})=\frac{\partial}{\partial z},~J\frac{\partial}{\partial z}=-(\frac{\partial}{\partial y}+x\frac{\partial}{\partial z}).\]
	The vector fields
    $$e_1=\frac{1}{2}\left(\frac{\partial}{\partial t}-i \frac{\partial}{\partial x}\right)\mbox{ and } e_2=\frac{1}{2}\left(\left(\frac{\partial}{\partial y}+x \frac{\partial}{\partial z}\right)-i\frac{\partial}{\partial z}\right),$$
    span $T^{1,0} \Omega$ at each point. It is clear that 
    \[div(e_1)=div(e_2)=0,\]
    and
    \begin{eqnarray*}
        [e_1, e_2]=
        [e_1, \bar e_2]=-[\bar e_1,e_2]=-[\bar e_1,\bar e_2]=\frac{i}{2}(e_2-\bar e_2).
    \end{eqnarray*}
\end{exa}

\subsection{Formula on Chern connection}
This subsection is about the $(1,1)$-regularity on the almost complex manifolds, which is related to \cite[Appendix C.3]{TWZZ}
First, we recall the calculation of the Chern connection on almost Hermitian manifolds.
Let $\nabla$ be the metric and almost complex preserving connection with vanishing $(1,1)$-torsion.   Let $\{Z_i\}$ be a local frame on $T^{1,0}(X)$. We set 
\[[Z_i,Z_j]=\sum_kB^k_{ij}Z_k+B^{\bar k}_{ij}Z_{\bar k},\]
\[[Z_i,Z_{\bar j}]=\sum_k B^k_{i\bar j}Z_k+B^{\bar k}_{i\bar j}Z_{\bar k},\]
where $B^{\bar k}_{ij}=N^{\bar k}_{ij}$. Since the $(1,1)$-torsion vanishes, i.e. 
$\Gamma^{\bar k}_{i\bar j}=B^{\bar k}_{i\bar j}$ and $\Gamma^{\bar k}_{\bar ji}=0$, it derives 
\begin{eqnarray*}
\nabla_{Z_i}Z_j&=&\sum_kg
(\nabla_{Z_i}Z_j,Z_{\bar k})Z_k\\
&=&\sum_kZ_i(g(Z_j,Z_{\bar k}))Z_k-g(Z_j,\nabla_{Z_i}Z_{\bar k})Z_k\\
&=&\sum_kZ_i(g(Z_j,Z_{\bar k}))Z_k-\sum_{k,l}g(Z_j,B^{\bar l}_{i\bar k}Z_{\bar l})Z_k.
\end{eqnarray*}
Hence, we obtain 
\begin{equation}
\Gamma^k_{ij}=g^{k\bar l}Z_i(g_{j\bar l})-g^{k\bar l}g_{j\bar r}B^{\bar r}_{i\bar l}.\label{eqn-chern-connection-structure}
\end{equation}
So, the authors can not follow how to get the formula \cite[Appendix (C.1)]{TWZZ}, since there is an  extra term in \eqref{eqn-chern-connection-structure}. As the last example in the previous subsection and by the same reason, the authors also can not follow how to derive \cite[Appendix Proposition C.4.]{TWZZ}, which is as the same as in the complex case(cf. \cite[Proposition 2.9]{Dem92}).

\section{The proof}

Evidently, the subspaces of smooth functions and forms are dense in the corresponding $L^2$-space. 
%Utilizing the model in Theorem \ref{lemma-partial-bar}, we set \[ H_1=L^2(M,\mathbb R)_0, ~H_2=L^2(\mathcal{A}^1_{\mbr})\mbox{ and } H_3=L^2(\mathcal{A}^-_J).\]It suffices to show that the complex\begin{eqnarray}H_1 \xrightarrow{\mathcal W_d} H_2    \xrightarrow{d_J^{-}}H_3,\label{complex-exact-1}\end{eqnarray}is exact. 
Before proceeding, %we review an elementary lemma. Readers who are familiar with can skip to the proof part.
%\begin{lemma}\label{lemma-sympl-vector-Hodge}Let $(V,\omega)$ be a $2n$-dimensional symplectic vector space, and $J$ be a compatible (almost)complex structure on $(V,\omega)$. Then, for any $\alpha\in V^*$, the identity \[*(J\alpha)=\alpha\wedge \frac{\omega^{n-1}}{(n-1)!},\]where $*$ is the Hodge operator with respect to the metric $\omega(-,J-)$. \end{lemma}\begin{pf}By the linear model of the previous section, we can find a unitary basis $\{e_1,e_2,\cdots e_{2n}\}$ of $V$ such that \[\omega=e^1\wedge e^2+e^3\wedge e^4+\cdots +e^{2n-1}\wedge e^{2n},\]and $J(e_{2i-1})=e_{2i},~J(e_{2i})=-e_{2i-1}$,  where $e^i$ is dual to $e_i$ for each $1\leq i\leq n$. It suffices to check the identity of the cases $\alpha=e^{2j-1}$ and $\alpha=e^{2j}$.The formulas $*J(e^{2j-1})=*e^{2j}=-e^1\wedge\cdots\widehat{e^{2j}}\wedge e^{2n}= e^{2j-1}\wedge\frac{\omega^{n-1}}{(n-1)!}$ and$*J(e^{2j})=-*e^{2j-1}=e^1\wedge\cdots\widehat{e^{2j-1}}\wedge e^{2n}=e^{2j}\wedge \frac{\omega^{n-1}}{(n-1)!}$ yield the identity of the lemma. \end{pf}
we review the basic property of the Atiyah-Hitchint-Singer operator \cite{AHS}, 
$$d^+\oplus d^*:\mcA^{1}_{\mbr}\to \mcA^+_\mbr\oplus \mcA^0_\mbr,$$ of a compact oriented $4$ manifold $(M,g)$ without boundary.
Since $d^+\oplus d^*$ is a first-order elliptic (Fredholm)operator, one has the regularity 
\begin{eqnarray}
    \|a\|_{L^2_1}\leq C(\|(d^+\oplus d^*)(a)\|_{L^2}+\|a_h\|_{L^2}),\label{formu-AHS-operator}
\end{eqnarray}
where $a_h$ denotes the Hodge-Laplacian harmonic component of $a$ and $C$ is some constant depending only on $(M,g)$. 
Here, the fact $\ker(d^+)=\ker (d)$ is provided by 
\[0=\int_Mda\wedge da=\int_{M}d^+a\wedge d^+a-d^-a\wedge d^-a=\|d^+a\|^2_{L^2}-\|d^-a\|^2_{L^2}.\]
Next, recall that for any Hermitian  almost complex $4$ manifold $(M,J,F)$, one has the decompositions, c.f., \cite[Lemma 2.1.57]{Don86},
\[\Lambda^+=\Lambda^{2,0}\oplus\Lambda^{0,2}\oplus\mathbb C\langle F\rangle,~\Lambda^-=\ker(\Lambda_F)\cap\Lambda^{1,1},\]
where  $\Lambda^+$ and $\Lambda^-$ are the $1$ and $-1$-eigenspaces of the Hodge-$*$ operator  with respect to the metric $F(-,J-)$ respectively. 
So, we write 
\begin{eqnarray}
    d^+a=d^-_Ja+\Lambda_Fd^+_J\cdot F.\label{formu-d+}
\end{eqnarray}
Recall that $\Lambda_F$ denotes the contraction with the form $F$.

Now, we give the proof. 

\begin{pf}\textbf{ of Theorem \ref{thm-1}}
 %   Recall the Hodge-decomposition gives the estimate, c.f. \cite[Lemma 5.1]{Weih},     \begin{eqnarray}      \|a\|_{L^2} \leq  \|a\|_{L^2_{1}}\leq C(\|da+d^*a)\|_{L^2}+\|a_{h}\|_{L^2}),\label{est-regulairty}    \end{eqnarray}    where $C$ is a constant depending only on $(M,J,F)$ and $a_h$ denotes the harmonic part of $a\in \mathcal A^1_\mbr$. 
 The idea is to construct a suitable model in Theorem \ref{lemma-partial-bar}.
 Assume that $\psi=d^+_Ja$ is non-zero with $d^-_Ja=0$. This is equivalent to say that $a\in(\ker d)^{\perp}$ by the Hodge decomposition,
 \[\mcA^1_\mbr=im(d)\oplus Harm\oplus im(d^*),\]
 where $Harm$ denotes the space of harmonic $1$-forms.
 It suffices to show that the complex, 
    \begin{eqnarray}
H_1 \xrightarrow{\tmw } H'_2    \xrightarrow{d_J^{-}}H_3,\label{complex-exact-2}
\end{eqnarray}
is exact, where 
\[ H_1= L^2(M,\mathbb R)_0 ,~ H_2=\{d^*b\in L^2\mcA^1_\mbr\mid b\in  \mcA^2_\mbr\} \mbox{ and }H'_3=L^2(\mathcal{A}^-_J) .\]
Recall that the form $F$ with point-wise norm squar $2$ with respect to the metric $F(-,J-)$.
%By the Hodge decomposition \[\mcA^1_\mbr=im(d)\oplus Harm\oplus im(d^*),\]we choose $V=Harm^\perp$ in the model of Theorem \ref{lemma-partial-bar}. Here $Harm$ denotes the space of harmonic $1$-forms.
    The formula  \eqref{formu-tW-adjoint} derives 
    \[\|\tmw^*(a)\|^2_{L^2}=\int_M |\Lambda d^+_Ja|^2\frac{F^2}{2}.\]
    By the formula \eqref{formu-d+} and the estimate \eqref{formu-AHS-operator}, we obtain 
    \begin{eqnarray}
        \|a\|_{L^2}&\leq& \|a\|_{L^2_1}\leq C(\|d^+a\|_{L^2}+\|d^*a\|_{L^2}+\|a_h\|_{L^2})\nonumber\\
        &=&C(\|d^-_Ja\|_{L^2}+\|\Lambda_Fd^+_Ja\cdot F\|_{L^2}+\|d^*a\|_{L^2}+\|a_h\|_{L^2})\nonumber\\
        &\leq &2C(\|d^-_Ja\|_{L^2}+\|\Lambda_Fd^+_Ja\|_{L^2}+\|d^*a\|_{L^2}+\|a_h\|_{L^2})\nonumber\\
        &=&2C(\|d^-_Ja\|_{L^2}+\|\tmw^*a\|_{L^2}+\|d^*a\|_{L^2}+\|a_h\|_{L^2}).\label{esti-wd-1}
    \end{eqnarray}
    When $a\in H_2$, the estimate \eqref{esti-wd-1} gives 
    \begin{eqnarray}
        \|a\|_{L^2}\leq C_1(\|d^-_Ja\|_{L^2}+\|\tmw^*a\|_{L^2}+\|d^*a\|_{L^2}),\label{est-wd-2}
    \end{eqnarray}
    where $C_1$ is a constant depending only on $(M,J,F)$. 
    Thus, \eqref{complex-exact-2} is exact by Theorem \ref{lemma-partial-bar}.
    
\begin{rmk}
    For any $a\in \mcA^1_\mbr$, we write $a=df_a+a_h+d^*a_b$. %There is a unique $a':=a_h+d^*a_b$ such that $a'\in\ker(d^*)$. Hence,    \[\|a'\|_{L^2_1}\leq C(\|d^+a'\|_{L^2}+\|a_h\|_{L^2}).\]
    %It is   clear  that  there is no such  a constant $C$ to hold  $\|a\|_{L^2}\leq C\|a'\|_{L^2_1}$. 
    In general, there is non constant $C$ to hold the inequality
   \[\|a\|_{L^2}
   \leq C({\|d^-_Ja\|_{L^2}+\|\mw^*a\|_{L^2}}),\]
   for a random $1$-form $a$.
   For example, on  acompact almost K\"ahler $4$ manifold, one can choose $a$ to be a harmonic $1$-form or $df$ for some non-constant $f$ to get a contradiction. 
   Hence, it is natural to consider  the Atiyah-Hitchint-Singer operator on the $d^*$-exact part to show the global exactness.  
  % The author also notices that in Wang Wang and Zhu use the operator $d+d^*$ on the whole space $L^2(\mcA^1_\mbr)$ to show the exactness
\end{rmk}

\end{pf}

%\begin{pf}\textbf{ of Theorem \ref{thm-2}}    Recall that for almost  K\"ahler $4$ manifold, we have $\mw^*(a)=\Lambda_Fda$.We argue by contradiction. Assume that  there is a constant $C$ to hold \begin{eqnarray}    \|a\|_{L^2}\leq C(\|d^-_Ja\|_{L^2}+\|\mw^*a\|_{L^2},\label{est-dJ-W}\end{eqnarray}for any $a\in \mcA^1_\mbr$.  Choose $a=df$ for some non-constant function $f$ or $a$ is a non-zero harmonic $1$-form. Put $a$ into the above inequality \eqref{est-dJ-W},the left-hand side $\|a\|_{L^2}>0$, while the right-hand side is vanishing. This contradicts to the assumption.\end{pf}

 %%%%%%%%%%%%%

College of Mathematics and Statistics, Chongqing University,Huxi Campus, Chongqing, 401331, P. R. ChinaChongqing Key Laboratory of Analytic Mathematics and Applications, Chongqing University, Huxi Campus, Chongqing, 401331, P. R.China 

E-mails: dexielin@cqu.edu.cn (or lindexie@126.com),

\quad  \quad \quad\quad  zhouhyu@cqu.edu.cn


\begin{thebibliography}{99}

\bibitem{AHS}
Atiyah,
M. F. ; Hitchint, N. J.; Singer,  I. M.: Self-duality in four-dimensional Riemannian geometry, Proc. R. Soc. Lond. A. 362, (1978),425-461.  

    \bibitem{CSZ}
    Chen, Zhihua; Siu, Yuntong; Zhong, Jiaqing: Introduction to Several Complex Variable,(in Chinese), High Educational Press.(2012) ISBN: 9787040362688

    \bibitem{Dem92}
Demailly, Jean-Pierre, Regularization of closed positive currents and intersection theory, J.
Alg. Geom., 1 (1992), 361-409.

\bibitem{DLiZh10}
Draghici, Tedi; Li, Tian-Jun and Zhang, Weiyi: Symplectic forms and cohomology decomposition of almost complex four-manifolds, Int. Math. Res. Not. IMRN (2010), no. 1, 1-17. MR 2576281

\bibitem{Don86}
Donaldson, Simon Kirwan; Kronheimer, Peter Benedict: The geometry of four manifolds, Oxford Math. Monogr.
Oxford Sci. Publ.
The Clarendon Press, Oxford University Press, New York, 1990. x+440 pp.
ISBN:0-19-853553-8
 
\bibitem{Don06}
Donaldson, S. K.: Two forms on four manifolds and elliptic equations, Nankai Tracts
Math., 11, Inspired by S. S. Chern, 153-172, World Sci. Publ., Hackensack, N.J.,
2006.


\bibitem{HV15}
Harvey, R.;  Lawson, Jr.H. B.: Potential theory on almost complex manifolds, Ann.
Inst. Fourier, 65 (2015), 171-210.

\bibitem{HLP16}
Harvey, F.R., Lawson, H.B. ,  Pliś, S.: Smooth approximation of plurisubharmonic functions on almost complex manifolds. Math. Ann. \textbf{366} (2016), 929-940

\bibitem{Hor65}
H\"ormander, L.: $L^2$ estimates and existence theorems for the $\bar\partial$-operator, Acta Math.,
113(1965), 89-152. 

\bibitem{Hor90}
H\"ormander, L.: An introduction to complex analysis in several variables, third edition
(revised), D. Van Nostrand Co., Inc., Princeton, N.J.-Toronto, Ont.-London, 1990.

\bibitem{Lejmi06}
Lejmi, M.: Strictly nearly K\"ahler 6-manifolds are not compatible with symplectic forms,
C. R. Math. Acad. Sci. Paris 343 (2006), no. 11-12, 759-762.



%\bibitem{Morrey66} C.B. Morrey,  Multiple Integrals in the Calculus of Variations , Grundlehren dermathematischen Wissenschaften, vol. 130, Springer-Verlag, Berlin, 1966.


%\bibitem{MW52} Motzkin, I. and Wasow, W. On the approximation of linear elliptic differential equations by difference equations with positive coefficients, J. Math. Phys. 31 (1952), 253-259.

%\bibitem{NN}A. Newlander, L.Nireburg, Complex analytic coordinates in almsot complex manifolds,Ann. of Math. \textbf{65} (1957), pp. 391-404.



%\bibitem{Silva} A. Silva, $\partial\bar\partial$-closed positive currents and special metrics on compact complex manifolds, Complex Analysisand Geometry, Lecture notes in Pure and Appl. Math. 173 (Dekker, New-York,1996) 377-441.

\bibitem{TWZZ}
Tan, Qiang; Wang,  Hongyu; Zhou, Jiuru and Zhu, Peng: On tamed almost complex fourmanifolds, Peking Math. J. 5 (2022), no. 1, 37-152. MR 4389489

%\bibitem{Usher07}M. Usher, Standard surfaces and nodal curves in symplectic 4-manifolds, J. DifferentialGeom. 77 (2007), no. 2, 23-290.

%\bibitem{Weih} Weiharm; Uhlenbeck compactness, EMS

\bibitem{Wein07}
Weinkove, Ben: The Calabi-Yau equation on almost K\"ahler four-manifolds, J. Diff.
Geom., 76 (2007), 317-349.

%%%%%%%%%%
\end{thebibliography}
\end{document}